\documentclass[11pt]{article}
\usepackage[a4paper,margin=28mm]{geometry}
\usepackage{amsmath,amssymb,amsthm}
\usepackage[T1]{fontenc}
\usepackage{lmodern}
\usepackage{enumitem}
\usepackage[hidelinks]{hyperref}
\newtheorem{theorem}{Theorem}[section]
\newtheorem{lemma}[theorem]{Lemma}
\newtheorem{proposition}[theorem]{Proposition}
\newtheorem{conjecture}[theorem]{Conjecture}
\newtheorem{corollary}[theorem]{Corollary}
\theoremstyle{remark}

\newcommand{\inn}{\operatorname{in}}
\newcommand{\out}{\operatorname{out}}
\newcommand{\defic}{\operatorname{def}}
\title{Strong complete immersions in graphs with\protect\\ independence number two}
\author{Zijian Deng$^{1}$ \and Chuanshu Wu$^{2}$\thanks{Corresponding author.\newline E-mail addresses: \mbox{\href{mailto:zj1329205716@163.com}{\texttt{zj1329205716@163.com}} (Z. Deng)}, \mbox{\href{mailto:wuchuanshu2021@163.com}{\texttt{wuchuanshu2021@163.com}} (C. Wu)}.}}
\date{\small $^{1}$School of Mathematics and Statistics, Lingnan Normal University,\\
 Zhanjiang, Guangdong 524048, China\\[4pt]
 $^{2}$Center for Discrete Mathematics, Fuzhou University,\\
 Fuzhou, 350108, Fujian, China}
\begin{document}
\maketitle
\begin{abstract}
We prove that every graph $G$ with independence number at most two contains a strong immersion of $K_{\chi(G)}$. This strengthens Dahlke's theorem for weak immersions and establishes the Lescure--Meyniel conjecture for this class of graphs.
\end{abstract}

\section{Introduction}
All graphs considered in this paper are finite and simple. An \emph{immersion} of a graph $F$ in a graph $G$ consists of an injection $\varphi:V(F)\to V(G)$ and pairwise edge-disjoint paths $P_{uv}$ joining $\varphi(u)$ and $\varphi(v)$ for $uv\in E(F)$. The vertices in $\varphi(V(F))$ are its \emph{branch vertices}. The immersion is \emph{strong} if
\[
V(P_{uv})\cap\varphi(V(F))=\{\varphi(u),\varphi(v)\}
\qquad(uv\in E(F)).
\]
We write $G\succeq_s F$ if $G$ contains a strong immersion of $F$.

Lescure and Meyniel proposed the following conjecture.
\begin{conjecture}[Lescure and Meyniel~{\cite{LM}}]\label{conj:LM}
Every graph $G$ contains a strong immersion of $K_{\chi(G)}$.
\end{conjecture}
Abu-Khzam and Langston~\cite{AKL} proposed the corresponding conjecture for weak immersions. We prove Conjecture~\ref{conj:LM} for graphs with independence number at most two.
\begin{theorem}\label{thm:main}
If $\alpha(G)\leq2$, then $G\succeq_s K_{\chi(G)}$.
\end{theorem}
For an $n$-vertex graph with $\alpha(G)\leq2$, we have $\chi(G)\geq\lceil n/2\rceil$. Thus Theorem~\ref{thm:main} also gives $G\succeq_s K_{\lceil n/2\rceil}$. Vergara~\cite{V} proved $G\succeq_s K_{\lceil n/3\rceil}$, and Gauthier, Le and Wollan~\cite{GLW} improved this to $G\succeq_s K_{2\lfloor n/5\rfloor}$. Botler et al.~\cite{B} proved $G\succeq_s K_{\chi(G)}$ when $n\geq11$ and $\Delta(G)<19n/29-1$, or when $G$ is covered by three cliques and $\Delta(G)<2n/3-1$. Deng, Chang and Liu~\cite{DCL} proved the gem-free case. Quiroz~\cite{Q} proved that every $F$-free graph $G$ with $\alpha(G)\leq2$ contains an immersion of $K_{\chi(G)}$, provided $|V(F)|\leq4$ and $\alpha(F)\leq2$.

Botler et al.~\cite{BJL} proved that every $n$-vertex graph with independence number two contains every complete bipartite graph on $\lceil n/2\rceil$ vertices as an immersion. Chen and Deng~\cite{CD} gave a simpler proof. For $\alpha(G)\geq3$, Bustamante et al.~\cite{BQSZ} obtained lower bounds on the order of a clique immersion in terms of $|V(G)|$ and $\alpha(G)$.

Dahlke~\cite{D} proved that every graph with $\alpha(G)=2$ contains a weak immersion of $K_{\chi(G)}$. Theorem~\ref{thm:main} strengthens this result to strong immersions. We use the bridge construction and the $L/M/R$ decomposition of~\cite{D}. By choosing the matching and the available paths suitably, we eliminate the $RR$ pieces, whose paths may contain other branch vertices.

Berm\'udez, Netto and Quiroz~\cite{N} proved that every $n$-vertex graph with $\alpha(G)=2$ contains $K_{\lceil n/2\rceil}$ as an odd minor or as a totally odd strong immersion; the immersion alternative holds under additional hypotheses. Echeverr\'ia and McDonald~\cite{E} proved $\chi(G)\leq\lceil3(t-1)/2\rceil$ for graphs in this class with no strong odd $K_t$ immersion. We impose no parity restriction on the immersion paths.

Our decomposition lemma uses Tan's stable partitions and Fleiner's predecessor-counting argument~\cite{T,F}. Section~\ref{sec:extension} proves an extension theorem, from which Theorem~\ref{thm:main} follows by induction in Section~\ref{sec:main}.

\subsection*{Notation and preliminaries}
Write $[t]=\{1,\ldots,t\}$, $xy\in G$ for $xy\in E(G)$, $d_G(x,S)=|N_G(x)\cap S|$, and $\nu(G)$ for the matching number. For disjoint sets $X,Y$, let $E_G(X,Y)$ be the edges between them. The join of disjoint graphs $G_1,G_2$ is denoted by $G_1\vee G_2$. A graph $J$ is \emph{factor-critical} if $J-x$ has a perfect matching for every $x\in V(J)$. We allow the empty $K_0$ immersion.

\begin{theorem}[Gallai--Edmonds; see {\cite[Theorem 5]{W}}]\label{thm:GE}
Let $D$ be the vertices missed by some maximum matching of $J$, and put
\[
A=N_J(D)\setminus D,\qquad C=V(J)\setminus(D\cup A).
\]
Then $J[C]$ has a perfect matching, every component of $J[D]$ is factor-critical, and $E_J(C,D)=\varnothing$. Every maximum matching consists of a perfect matching on $C$, near-perfect matchings in the components of $J[D]$, and edges matching $A$ into distinct such components. In particular,
\[
\defic(J):=|V(J)|-2\nu(J)=c(J[D])-|A|.
\]
\end{theorem}

The following local-order form follows from the increasing and decreasing colour orders in the proof of Borodin, Kostochka and Woodall~\cite[Theorem 3, p.~189]{BKW}.
\begin{theorem}[Borodin, Kostochka and Woodall]\label{thm:BKW}
A simple bipartite graph $B$ admits strict orders on its incident edges such that
\[
p_x(xy)+p_y(xy)\leq\max\{d_B(x),d_B(y)\}-1
\qquad(xy\in E(B)),
\]
where $p_x(e)$ is the number of edges preceding $e$ at $x$.
\end{theorem}
Since $B$ is simple, the seeing count in~\cite{BKW} equals the predecessor sum above.

\begin{theorem}[Tan~{\cite{T}}]\label{thm:SP}
In a graph $J$ with strict preferences on the neighbours of each vertex and self-partnership ranked last, there is a permutation $\pi$ such that $x\pi(x)\in E(J)$ whenever $\pi(x)\ne x$,
\[
\pi(x)\succeq_x\pi^{-1}(x),
\]
and no edge $xy$ satisfies both $y\succ_x\pi^{-1}(x)$ and $x\succ_y\pi^{-1}(y)$.
\end{theorem}
See also~\cite[Theorem 2]{BCK}. Replacing an even cycle of $\pi$ by alternating pairs gives each vertex a partner at least as preferred as its predecessor. Thus we may assume that every cycle of $\pi$ has length one, two, or odd length at least three.

\section{Local orders and decompositions}
Let $Q$ be a simple graph and $U$ a finite colour set. For each $i\in V(Q)$, put
\[
U=L_i\sqcup M_i\sqcup R_i,\qquad A_i=L_i\cup M_i,\qquad r_i=|A_i|.
\]
The following pieces, each assigned a colour $w\in U$, are adapted from~\cite[Lemma 2.7]{D}:
\begin{enumerate}[label=(\roman*),topsep=3pt,itemsep=2pt]
\item an $F$-edge $ij$ with owner $i$ and $w\in M_i$, using the pair $(i,w)$;
\item an $LL$- or $LR$-edge $ij$, using $(i,w),(j,w)$;
\item an odd cycle with all vertices labelled $L$, using $(i,w)$ at each vertex;
\item a path $a-b-c$ labelled $L,L,R$, using $(a,w),(b,w),(c,w)$.
\end{enumerate}
Here a vertex has label $L$, $M$ or $R$ according to the set containing $w$. A \emph{compatible decomposition} partitions $E(Q)$ into these pieces and uses each pair $(i,w)$ at most once.

\begin{lemma}\label{lem:pieces}
If $Q$ has strict local edge orders satisfying
\begin{equation}\label{eq:ranks}
p_i(ij)+p_j(ij)<|A_i\cup A_j|\qquad(ij\in E(Q)),
\end{equation}
then $Q$ has a compatible decomposition.
\end{lemma}
\begin{proof}
Process the colours successively. For a colour $w$, let $Q_w$ consist of the remaining edges $ij$ with $w\in A_i\cup A_j$. Apply Theorem~\ref{thm:SP} to $Q_w$, preferring earlier incident edges.

Assign each 2-cycle an $F$-edge if it has an $M$ endpoint, and an $LL$- or $LR$-edge otherwise. Retain every all-$L$ odd cycle. Let $x_0x_1\cdots x_{2t}x_0$ be any other odd cycle, in successor order.

If $x_0$ has label $M$, select the $F$-edge $x_0x_1$ owned by $x_0$ and the pairs
\[
x_1x_2,\ x_3x_4,\ldots,\ x_{2t-1}x_{2t}.
\]
If the cycle has no $M$ vertex, its $R$ vertices are independent. Hence, after relabelling, $x_0,x_1,x_2$ have labels $L,L,R$. Select this path and the pairs $x_3x_4,\ldots,x_{2t-1}x_{2t}$. Assign the pairs as above and delete all selected edges. The selected pieces use distinct pairs $(i,w)$.

Every unselected edge $e$ of $Q_w$ has a selected edge preceding it at one endpoint. If one endpoint of $e$ is a singleton, then, since self-partnership is ranked last, Theorem~\ref{thm:SP} forces the other endpoint to be non-singleton, with both its predecessor and successor edges preceding $e$. Otherwise, if $e$ is outside the permutation cycles, the same theorem gives an endpoint at which both cycle edges precede $e$. In either case, at least one of these cycle edges is selected. Finally, if $e=x_ix_{i+1}$ is an omitted cycle edge, then $x_{i+1}x_{i+2}$ is selected and precedes $e$ at $x_{i+1}$ (subscripts modulo $2t+1$).

If $ij$ remains after all colours, each $w\in A_i\cup A_j$ therefore gives a distinct deleted predecessor at $i$ or $j$. Thus
\[
p_i(ij)+p_j(ij)\geq|A_i\cup A_j|,
\]
contrary to~\eqref{eq:ranks}.
\end{proof}

Assume henceforth in this section that $r_i\geq d_Q(i)$. Call $i$ \emph{tight} if $r_i=d_Q(i)$ and \emph{slack} otherwise. An orientation is \emph{nearly balanced} if $|\out(i)-\inn(i)|\leq1$ for every $i$. Such an orientation is obtained by pairing the odd vertices in each component, adding the corresponding edges, orienting an Euler tour, and deleting the added edges. An odd vertex is a \emph{source} if $\out(i)=\inn(i)+1$ and a \emph{sink} otherwise.

\begin{lemma}\label{lem:middle}
Let $M$ be a matching whose endpoints have odd $Q$-degree. Orient $Q-M$ nearly balanced and put $\varepsilon_i=\mathbf1_{i\in V(M)}$. There are local orders with incoming edges first, the matching edge next, and outgoing edges last, such that every arc $i\to j$ satisfies
\begin{equation}\label{eq:middle}
p_i(ij)+p_j(ij)\leq\inn(i)+\varepsilon_i+
\max\{\out(i),\inn(j)\}-1.
\end{equation}
The right-hand side is at most $\max\{d_Q(i),d_Q(j)\}-1$, except when $d_Q(i)=d_Q(j)$ is odd, $j$ is an unmatched sink, and $i$ is matched or an unmatched sink. In this case it equals $d_Q(i)$. For $ij\in M$, the predecessor sum is at most $\max\{d_Q(i),d_Q(j)\}-1$.
\end{lemma}
\begin{proof}
Split each vertex into an outgoing copy and an incoming copy, and apply Theorem~\ref{thm:BKW} to the resulting simple bipartite graph. For $i\to j$, the incoming edges and the matching edge at $i$ add $\inn(i)+\varepsilon_i$ predecessors, giving~\eqref{eq:middle}.

If $\out(i)\geq\inn(j)$, the bound is $d_Q(i)-1$. Otherwise it exceeds both endpoint bounds only if
\[
\inn(j)>\out(i),\qquad
\inn(i)+\varepsilon_i>\out(j)+\varepsilon_j.
\]
Set $\delta_i=\inn(i)+\varepsilon_i-\out(i)$ and $\eta_i=\inn(i)-\out(i)-\varepsilon_i$. Their values are
\[
\begin{array}{c|rrrr}
i&\text{matched}&\text{unmatched sink}&\text{even}&\text{unmatched source}\\\hline
\delta_i&1&1&0&-1\\
\eta_i&-1&1&0&-1
\end{array}
\]
The two strict inequalities force $\delta_i=\eta_j=1$ and $\out(i)=\out(j)$, giving the stated exception. For $ij\in M$,
\[
p_i(ij)+p_j(ij)=\frac{d_Q(i)-1}{2}+\frac{d_Q(j)-1}{2}
\leq\max\{d_Q(i),d_Q(j)\}-1.
\]
\end{proof}
In the exceptional case, a slack endpoint gives $|A_i\cup A_j|\geq d_Q(i)+1$, so~\eqref{eq:ranks} still holds.

\begin{corollary}\label{cor:matching}
Suppose $Q$ has a matching covering all odd tight vertices except at most one, and every even-degree endpoint of the matching is slack. Then $Q$ has local orders satisfying~\eqref{eq:ranks}.
\end{corollary}
\begin{proof}
Add a pendant edge at each even matching endpoint. At such a vertex, $r_i\geq d_Q(i)+1$, so its unchanged set $A_i$ still has size at least its new degree. Give each new vertex a set $A_i$ of size at least two. All matching endpoints are now odd and all new vertices are slack.

Orient the graph minus the matching nearly balanced, making the unique unmatched odd tight vertex a source if it exists. Lemma~\ref{lem:middle} has no exceptional edge with two tight endpoints. Restrict the resulting orders to $Q$; the predecessor counts decrease and the original sets $A_i$ are unchanged.
\end{proof}

\begin{lemma}\label{lem:multipartite}
Let $O$ be the odd tight vertices of $Q$. If $Q[O]$ is complete multipartite, then $Q$ has local orders satisfying~\eqref{eq:ranks}.
\end{lemma}
\begin{proof}
The case $O=\varnothing$ follows from Corollary~\ref{cor:matching}. Suppose a part $A$ satisfies $|A|\geq|O\setminus A|$. Match $O\setminus A$ into $A$, orient the graph minus these edges nearly balanced, and restore them directed toward $A$. Their endpoints have even residual degree, so every vertex of $O\setminus A$ becomes a source. All tight odd sinks lie in the independent set $A$. Apply Lemma~\ref{lem:middle} with $M=\varnothing$.

Otherwise $Q[O]$ has a matching missing at most one vertex. For even $|O|$, repeatedly pair vertices in two largest parts; for odd $|O|$, first delete a vertex from a largest part. Corollary~\ref{cor:matching} applies.
\end{proof}

\section{An extension theorem}\label{sec:extension}
\begin{theorem}\label{thm:merge}
Let $H$ be a triangle-free graph of order $2k-1$, $G=\overline H$, and $d_H(v)=\Delta(H)=D$. Suppose $H-v$ has a perfect matching, with pairs $p_ic_i$ for $p_i\in P=N_H(v)$, $i\in[D]$. Put
\[
C=\{c_i:i\in[D]\},\quad Y=V(H)\setminus(\{v\}\cup P\cup C),\quad m=k-1-D.
\]
If $G[Y]\succeq_s K_m$, then $G\succeq_s K_k$.
\end{theorem}
Throughout this section, use the notation of Theorem~\ref{thm:merge}. Since $H$ is triangle-free, $P$ is independent and $P\cap C=\varnothing$. The matching on $Y$ has $m$ edges. Fix a strong $K_m$ immersion in $G[Y]$ with branch set $W$, and put $Z=Y\setminus W$. Thus $|W|=|Z|=m$.

\subsection{Resources and optimal selections}\label{subsec:resources}

For $i\in[D]$, define
\[
\begin{aligned}
L_i&=N_H(c_i)\cap W,& M_i&=N_H(p_i)\cap W,\\
A_i&=L_i\cup M_i,& R_i&=W\setminus A_i,\qquad r_i=|A_i|.
\end{aligned}
\]
Since $p_ic_i\in H$, we have $L_i\cap M_i=\varnothing$. The available \emph{index resources} and \emph{individual resources} at $i$ are, respectively,
\[
\begin{aligned}
\mathcal A_i&=\{j\in[D]\setminus\{i\}:c_ic_j,c_ip_j\in G\},\\
\mathcal U_i&=\{z\in Z:c_iz,p_iz\in G\}.
\end{aligned}
\]
Write $i\to j$ and $i\to z$ for these resources, and set $a_i=|\mathcal A_i|$, $u_i=|\mathcal U_i|$.

\begin{lemma}\label{lem:identity}
For each $i\in[D]$,
\begin{align}
a_i&=D-d_H(c_i)+d_H(c_i,Y),\notag\\
u_i&=m-d_H(c_i,Z)-d_H(p_i,Z),\notag\\
q_i:=a_i+u_i-r_i&=D-d_H(c_i)+m-d_H(p_i,Y)\geq0.
\label{eq:resources}
\end{align}
\end{lemma}
\begin{proof}
For $j\ne i$, triangle-freeness gives $|N_H(c_i)\cap\{p_j,c_j\}|\leq1$, and $j\in\mathcal A_i$ precisely when this intersection is empty. Since $p_i\in N_H(c_i)$ and $v\notin N_H(c_i)$,
\[
a_i=(D-1)-(d_H(c_i)-d_H(c_i,Y)-1).
\]
Also, $N_H(c_i)\cap N_H(p_i)=\varnothing$, yielding the formula for $u_i$. Now
\[
r_i=d_H(c_i,W)+d_H(p_i,W),\qquad
d_H(c_i)\leq D,\qquad d_H(p_i,Y)\leq m.
\]
The last inequality follows because $N_H(p_i)$ meets each matching edge of $H[Y]$ at most once. Substitution proves~\eqref{eq:resources}.
\end{proof}

Choose a set $S_i\subseteq\mathcal A_i\sqcup\mathcal U_i$ with $|S_i|=r_i$ for each $i$. Here indices and vertices of $Z$ are treated as distinct labels. Define $Q$ on $[D]$ by
\[
ij\in E(Q)\quad\Longleftrightarrow\quad j\in S_i\text{ and }i\in S_j.
\]
Then $d_Q(i)\leq r_i$. We vary the perfect matching of $H[P,C]$ and the sets $S_i$, keeping $P,C,Y,W,Z$ and the labels $p_i$ fixed.

Put
\[
B=\{i:L_i=\varnothing\},\qquad I=[D]\setminus B,\qquad R=Q[I].
\]
Assign each edge of $Q$ incident with $B$ to an endpoint $i\in B$. Since $d_Q(i)\leq r_i=|M_i|$, distinct colours in $M_i$ give $F$-edges using no pair $(j,w)$ with $j\in I$. It remains to decompose $R$.

For $i\in I$, let $\mathcal F_i=(\mathcal A_i\cap B)\sqcup\mathcal U_i$. Call the selection \emph{normalised} if
\[
|S_i\cap\mathcal F_i|=\min\{r_i,|\mathcal F_i|\}\qquad(i\in I).
\]
Replacing an element of $S_i\cap I$ by an unused element of $\mathcal F_i$ cannot increase $|E(R)|$. Choose a matching and a normalised selection minimising
\begin{equation}\label{eq:potential}
(\mu,|T|)\quad\text{lexicographically},\qquad
\mu=|E(R)|,\quad T=\{i\in I:d_R(i)=r_i\}.
\end{equation}
If $i\in T$, then $r_i>0$ and every element of $S_i$ is a mutual index in $I$. Hence
\[
\mathcal F_i=\varnothing,\qquad
j\in\mathcal A_i\ \Longrightarrow\ i\in S_j.
\]
Indeed, a violation of the second implication permits replacing a mutual index by $j$, decreasing $\mu$. Consequently, for distinct $i,j\in T$,
\[
j\in\mathcal A_i\quad\Longleftrightarrow\quad i\in\mathcal A_j
\quad\Longleftrightarrow\quad ij\in E(R).
\]

\subsection{Exchanging matched vertices}
\begin{lemma}\label{lem:transport}
Let $a,b\in T$ satisfy $c_ac_b\in G$, $b\notin\mathcal A_a$ and $a\notin\mathcal A_b$. Then $p_ac_a,p_bc_b$ can be replaced by $p_ac_b,p_bc_a$, with a normalised selection for which $\mu$ does not increase. Moreover, for $c\in\{c_a,c_b\}$ and $j\notin\{a,b\}$,
\[
cc_j,cp_j\in G\quad\Longrightarrow\quad c_jp_a,c_jp_b\in G.
\]
\end{lemma}
\begin{proof}
Put $P_0=\{p_a,p_b\}$ and $C_0=\{c_a,c_b\}$. The assumptions give
\[
H[P_0\cup C_0]=K_{2,2}\quad\text{with parts }P_0,C_0.
\]
Since $\mathcal U_a=\mathcal U_b=\varnothing$, every $z\in Z$ meets both matching edges in $H$. Triangle-freeness implies
\[
N_H(z)\cap(P_0\cup C_0)\in\{P_0,C_0\}.
\]
Thus, after the exchange, $\mathcal U'_a=\mathcal U'_b=\varnothing$. Resources into $B$ depend only on the $C$ endpoint and the unchanged pairs in $B$, so $\mathcal F'_a=\mathcal F'_b=\varnothing$. Furthermore,
\[
L'_a=L_b,\quad L'_b=L_a,\quad
M'_a=M_a,\quad M'_b=M_b,\quad r'_a+r'_b=r_a+r_b>0,
\]
and $r'_a,r'_b>0$. In particular, $I$ and $B$ are unchanged.

For $j\notin\{a,b\}$, transfer $S_j\cap\{a,b\}$ as follows. If $N_H(c_j)\cap P_0\ne\varnothing$, triangle-freeness gives $N_H(c_j)\cap C_0=\varnothing$; keep each selected index at its $P$ endpoint. Otherwise keep it with its $C$ endpoint. The selected resources remain distinct and available, and $|S'_j|=|S_j|$. For $j\in I\setminus\{a,b\}$, also $\mathcal F'_j=\mathcal F_j$ and $S'_j\cap\mathcal F_j=S_j\cap\mathcal F_j$.

By Lemma~\ref{lem:identity}, choose $|S'_a|=r'_a$ and $|S'_b|=r'_b$. Since $ab\notin E(R')$ and $R'-\{a,b\}=R-\{a,b\}$,
\[
\begin{aligned}
|E(R')|-|E(R)|
&=d_{R'}(a)+d_{R'}(b)-d_R(a)-d_R(b)\\
&\leq r'_a+r'_b-r_a-r_b=0.
\end{aligned}
\]

For the last assertion, suppose $cc_j,cp_j\in G$ but $c_jp_t\in H$ for some $t\in\{a,b\}$. Pair $c$ with $p_t$, exchanging the two pairs if necessary. The resource $t\to j$ remains available, but $j\to t$ is unavailable. Choose $j\in S'_t$, and complete $S'_a,S'_b$ by Lemma~\ref{lem:identity}. Then
\[
d_{R'}(a)+d_{R'}(b)\leq r'_a+r'_b-1=r_a+r_b-1,
\]
contrary to the minimality of $\mu$.
\end{proof}

\begin{lemma}\label{lem:anchor}
If $d_H(p_h,Y)=m$ for some $h$, then there are a strong $K_m$ immersion in $G[Y]$, a perfect matching of $H[P,C]$, and a normalised selection for which $Q$ has a compatible decomposition.
\end{lemma}
\begin{proof}
Set $W^*=N_H(p_h)\cap Y$. This is an independent transversal of the matching of $H[Y]$, so $G[W^*]=K_m$. Use this clique as the immersion, set $W=W^*$ and $Z=Y\setminus W$, redefine the resource sets, and minimise~\eqref{eq:potential} again. For every partner $c_h$ of $p_h$, triangle-freeness gives $L_h=\varnothing$; hence $h\in B$.

For $i\in T$, choose $w\in L_i$. Since $c_iw,p_hw\in H$, we have $c_ip_h\in G$. As $h\notin\mathcal A_i$, it follows that $c_ic_h\in H$. Thus
\[
\{c_i:i\in T\}\subseteq N_H(c_h)
\]
is independent in $H$. If $a,b\in T$ and $ab\notin E(R)$, Lemma~\ref{lem:transport} applies. For $j\in N_{R[T]}(a)\setminus N_{R[T]}(b)$, we have $c_ac_j,c_ap_j\in G$ but $c_jp_b\in H$, a contradiction. Therefore
\[
ab\notin E(R[T])\quad\Longrightarrow\quad
N_{R[T]}(a)=N_{R[T]}(b).
\]
Hence $R[T]$ is complete multipartite. The odd tight vertices of $R$ form a subset of $T$, so the subgraph they induce is also complete multipartite. Lemmas~\ref{lem:multipartite} and~\ref{lem:pieces} therefore give a compatible decomposition of $R$. Adding the $F$-edges assigned to $B$ gives a compatible decomposition of $Q$.
\end{proof}
By~\eqref{eq:resources}, $q_i=0$ implies $d_H(c_i)=D$ and $d_H(p_i,Y)=m$, so Lemma~\ref{lem:anchor} applies whenever some $q_i=0$.

\subsection{Auxiliary matchings}
Colour the pairs in $T$ by
\[
ij:\quad
\begin{cases}
\text{red},&ij\in E(R),\\
\text{blue},&c_ic_j\in H,\\
\text{white},&\text{otherwise}.
\end{cases}
\]
The blue graph is triangle-free, and every white pair satisfies Lemma~\ref{lem:transport}. There is no triangle with one red and two white edges: if $aj$ is red and $ab,bj$ are white, then $c_ac_j,c_ap_j\in G$ but $c_jp_b\in H$, contradicting that lemma.

\begin{lemma}\label{lem:threecolor}
On any induced vertex set, the red graph has at most two nonbipartite components. If it has two, these are its only components.
\end{lemma}
\begin{proof}
Let $K$ be a nonbipartite red component. For $x\notin V(K)$, put
\[
B_x=\{u\in V(K):xu\text{ is blue}\}.
\]
Each $B_x$ is a vertex cover of $K$, since an uncovered red edge would form a red--white--white triangle with $x$. If $xy$ is blue with $x,y\notin V(K)$, then $B_x\cap B_y=\varnothing$, and every edge of $K$ joins $B_x$ to $B_y$. This makes $K$ bipartite, a contradiction.

Thus there are no blue edges outside $K$. If the red graph outside $K$ had a component containing an edge $uv$ and another component containing a vertex $x$, then $uv,ux,vx$ would again be red, white, white. The assertion follows.
\end{proof}

Set
\[
O=\{i\in T:d_R(i)\text{ is odd}\},\qquad J=R[O],
\]
and apply Theorem~\ref{thm:GE} to $J$, with parts $D_0,A_0,C_0$. A nontrivial connected factor-critical graph is nonbipartite: a bipartition with sizes $a,b$ would require both $a-1=b$ and $a=b-1$. Lemma~\ref{lem:threecolor}, applied on $O\setminus A_0$, gives at most two nontrivial components of $J[D_0]$. If there are two, then
\[
C_0=\varnothing,\qquad c(J[D_0])=2,\qquad \defic(J)=2-|A_0|.
\]
\begin{proposition}\label{prop:ordinary}
If $\defic(J)\leq1$ or $J[D_0]$ has at most one nontrivial component, then $R$ has local orders satisfying~\eqref{eq:ranks}, and $Q$ has a compatible decomposition.
\end{proposition}
\begin{proof}
If $\defic(J)\leq1$, a maximum matching of $J$ and Corollary~\ref{cor:matching} give~\eqref{eq:ranks}. Otherwise $J[D_0]$ has at most one nontrivial component.

\smallskip\noindent\textbf{Case 1.} Every component of $J[D_0]$ is a singleton.

Let $M$ be the perfect matching on $C_0$ and $F$ the matching from $A_0$ to $D_0$ in a maximum matching of $J$. Orient $R-M-F$ nearly balanced, and restore $F$ from $A_0$ to $D_0$. The endpoints of $F$ have even degree in $R-M-F$, so the resulting orientation of $R-M$ is nearly balanced and every vertex of $A_0$ is a source.

Every unmatched tight odd sink lies in $D_0$, and
\[
E(J[D_0])=E_J(D_0,C_0)=\varnothing.
\]
Its tight odd neighbours therefore lie in $A_0$ and are sources. Lemma~\ref{lem:middle} has no exceptional edge with two tight endpoints.

\smallskip\noindent\textbf{Case 2.} $J[D_0]$ has exactly one nontrivial component $K$.

Choose $b\in V(K)$ and a maximum matching of $J$ missing $b$. By Theorem~\ref{thm:GE}, no matching edge from $A_0$ enters $K$. Let $F$ be its edges from $A_0$ to the singleton components of $D_0$, and let $M$ consist of its perfect matchings on $C_0$ and $K-b$.

Orient $R-M-F$ nearly balanced, making $b$ a source by reversing its component if necessary. Restore $F$ from $A_0$ to $D_0$. Every vertex of $A_0$ becomes a source, since its residual degree was even. Every unmatched tight odd sink lies in a singleton component of $D_0$ and has tight odd neighbours only in $A_0$. Again Lemma~\ref{lem:middle} has no exceptional edge with two tight endpoints.

In both cases, any exceptional edge has a slack endpoint, so~\eqref{eq:ranks} holds. Lemma~\ref{lem:pieces} gives a compatible decomposition of $R$; adding the $F$-edges assigned to $B$ gives one of $Q$.
\end{proof}

\begin{proposition}\label{prop:two}
Suppose the optimal choice satisfies $q_i>0$ for every $i\in[D]$. If $J[D_0]$ has two nontrivial components and $\defic(J)=2$, then a normalised selection can be chosen, with the same matching and branch set $W$, for which $Q$ has a compatible decomposition.
\end{proposition}
\begin{proof}
Let $K_1,K_2$ be the two nontrivial components. Then
\[
A_0=C_0=\varnothing,\qquad O=V(K_1)\sqcup V(K_2).
\]
Choose $bj\in E(K_1)$. Since $b\in T$ and $q_b>0$, we have $\mathcal F_b=\varnothing$ and $a_b-r_b=q_b>0$. Choose $x\in\mathcal A_b\setminus S_b$. By optimality, $b\in S_x$. Also $x\notin T$, since otherwise $bx\in E(R)$ would force $x\in S_b$.

Replace $S_b$ by $(S_b\setminus\{j\})\cup\{x\}$. The selection remains normalised, and
\[
R'=R-bj+bx,\qquad |E(R')|=|E(R)|,
\]
\[
d_{R'}(j)=r_j-1,\qquad d_{R'}(x)=d_R(x)+1.
\]
Since $j$ leaves $T$, the minimality of $|T|$ forces $x$ to enter $T$. Hence
\[
d_R(x)=r_x-1,\qquad T'=(T\setminus\{j\})\cup\{x\}.
\]
The vertex $j$ has even $R'$-degree and is slack. Fix $b_2\in V(K_2)$. For a factor-critical graph $K$ and $u\in V(K)$, let $M(K-u)$ denote a perfect matching of $K-u$.

\smallskip\noindent\textbf{Case 1.} $d_{R'}(x)$ is even.

The odd tight vertices of $R'$ are $O\setminus\{j\}$. The matching
\[
M'=M(K_1-j)\cup M(K_2-b_2)
\]
lies in $R'$ and covers all of them except $b_2$.

\smallskip\noindent\textbf{Case 2.} $d_{R'}(x)$ is odd.

The odd tight vertices of $R'$ are $(O\setminus\{j\})\cup\{x\}$. The matching
\[
M'=M(K_1-b)\cup M(K_2-b_2)\cup\{bx\}
\]
lies in $R'$: $x\notin O$, and none of its edges is $bj$. It covers all odd tight vertices except $b_2$, and its only even endpoint is the slack vertex $j$.

In each case, Corollary~\ref{cor:matching} and Lemma~\ref{lem:pieces} give a compatible decomposition of $R'$. Adding the $F$-edges assigned to $B$ gives a compatible decomposition of the final graph $Q$.
\end{proof}

\subsection{Realising the decomposition}
\begin{proposition}\label{prop:realisation}
Let $W$ be the branch set of a strong $K_m$ immersion in $G[Y]$. For any perfect matching of $H[P,C]$ and resource sets $S_i$ as above, if $Q$ has a compatible decomposition, then $G$ has a strong $K_k$ immersion with branch set $\{v\}\cup C\cup W$.
\end{proposition}
\begin{proof}
As in~\cite[Section 3]{D}, take
\[
\mathcal B=\{v\}\cup C\cup W
\]
as the branch set, retaining the strong $K_m$ immersion on $W$. Since $N_H(v)=P$, join $v$ directly to $C\cup W$. For $i<j$, use the edge $c_ic_j$ if it lies in $G$, and otherwise the path
\begin{equation}\label{eq:baseline}
c_i-p_j-p_i-c_j.
\end{equation}
The latter exists by triangle-freeness of $H$ and independence of $P$ in $H$. Also use every edge of $E_G(C,W)$. The remaining pairs are
\[
\{(c_i,w):i\in[D],\ w\in L_i\}.
\]
For $w\in L_i$, the available resources give paths
\begin{equation}\label{eq:long}
w-p_i-p_j-c_i\quad(j\in S_i\cap[D]),
\end{equation}
or $w-p_i-z-c_i$ for $z\in S_i\cap Z$. Here $wp_i\in G$ because $L_i\cap M_i=\varnothing$.

For a piece of colour $w$, make the following assignments:
\begin{enumerate}[label=(\roman*),topsep=3pt,itemsep=3pt]
\item For an $F$-edge $ij$ owned by $i$, delete $i\to j$. Since $w\in M_i$, the edge $c_iw$ has already been used.
\item For an $LL$-edge $ij$, assign $w-p_j-c_i$ and $w-p_i-c_j$, using $i\to j$ and $j\to i$.
\item For an $LR$-edge $ij$ with label $L$ at $i$, assign $w-p_j-c_i$, using $i\to j$.
\item For an all-$L$ odd cycle $i_0i_1\cdots i_{2t}i_0$, assign $w-p_{i_{s+1}}-c_{i_s}$ using $i_s\to i_{s+1}$ for each $s$, with subscripts modulo $2t+1$.
\item For an $LLR$ path $a-b-c$, assign $w-p_b-c_a$ and $w-p_c-c_b$, using $a\to b$ and $b\to c$.
\end{enumerate}
These short paths exist because their internal vertices have label $L$ or $R$. Let $f_i$ be the number of deleted resources at $i$, and $s_i$ the number used by short paths starting at a vertex of $L_i$ and ending at $c_i$. Compatibility gives
\[
f_i\leq|M_i|,\qquad s_i\leq|L_i|,\qquad
r_i-f_i-s_i\geq|L_i|-s_i.
\]
Assign the remaining $|L_i|-s_i$ pairs to distinct unused elements of $S_i$, using the corresponding long paths. Discard unused resources.

We check edge-disjointness by edge type.
\begin{enumerate}[label=(\roman*),topsep=3pt,itemsep=3pt]
\item Each edge $c_ip_j$ determines a unique resource $i\to j$, used at most once. Its occurrence in~\eqref{eq:baseline} requires $c_ic_j\in H$, whereas a resource path requires $c_ic_j\in G$.
\item An edge $p_ip_j$ can occur in two resource paths only if $ij\in E(Q)$. Each piece deletes or shortens at least one of $i\to j,j\to i$. Its occurrence in~\eqref{eq:baseline} is excluded by the same adjacency condition as in (i).
\item A short path using $wp_i$ uses the pair $(i,w)$, so there is at most one such path. If $w\in L_i$, the same piece assigns the pair $(c_i,w)$ a short path, excluding its long path. If $w\in R_i$, the edge $c_iw$ has already been used; if $w\in M_i$, then $wp_i\notin G$. Distinct long paths use distinct pairs $(i,w)$.
\item Each edge $p_iz$ or $c_iz$ determines a unique individual resource $i\to z$, used at most once. These edges belong to no other new path.
\item The old paths use only $E(G[Y])$, which no new path uses. The direct edges from $v$, within $C$, and between $C$ and $W$ occur in no other path.
\end{enumerate}
Thus the paths are pairwise edge-disjoint. Every new internal vertex lies in $P\cup Z$, and
\[
(P\cup Z)\cap\mathcal B=\varnothing,\qquad |\mathcal B|=1+D+m=k.
\]
Therefore $G\succeq_s K_k$.
\end{proof}

\begin{proof}[Proof of Theorem~\ref{thm:merge}]
If $m=0$, the direct edges from $v$ to $C$, the edges of $G[C]$, and the paths~\eqref{eq:baseline} give a strong $K_k$ immersion. Assume $m>0$ and choose a matching and a normalised selection minimising~\eqref{eq:potential}.

If some $q_i=0$, Lemma~\ref{lem:anchor} gives a choice of $W$, matching and selection for which $Q$ has a compatible decomposition. Otherwise $q_i>0$ for every $i$. By Lemma~\ref{lem:threecolor} and Theorem~\ref{thm:GE}, $J[D_0]$ has at most two nontrivial components. Proposition~\ref{prop:ordinary} applies if it has at most one, or if $\defic(J)\leq1$. In the remaining case, $\defic(J)=2-|A_0|>1$, so $A_0=\varnothing$ and $\defic(J)=2$; Proposition~\ref{prop:two} applies. Thus in every case a compatible decomposition of $Q$ exists. Proposition~\ref{prop:realisation} completes the proof.
\end{proof}

\section{Proof of the main theorem}\label{sec:main}
\begin{proof}[Proof of Theorem~\ref{thm:main}]
Let $H=\overline G$. Since $\alpha(G)\leq2$, the two-vertex colour classes form a matching in $H$, and
\begin{equation}\label{eq:chi}
\chi(G)=|V(G)|-\nu(H).
\end{equation}
We proceed by induction on $|V(G)|$. The assertion is immediate for $|V(G)|\leq1$. Passing to a vertex-critical induced subgraph with the same chromatic number, assume $G$ is vertex-critical. For every $v\in V(G)$,
\[
\chi(G-v)=\chi(G)-1\quad\Longrightarrow\quad\nu(H-v)=\nu(H).
\]
Hence $D=V(H)$ in Theorem~\ref{thm:GE}, so every component of $H$ is factor-critical.

If $H$ has components $H_1,\ldots,H_t$ with $t\geq2$, then
\[
G=\overline H_1\vee\cdots\vee\overline H_t,\qquad
\chi(G)=\sum_{i=1}^t\chi(\overline H_i).
\]
By induction, $\overline H_i\succeq_s K_{\chi(\overline H_i)}$. Joining their branch sets by the direct edges between different factors gives $G\succeq_s K_{\chi(G)}$.

Assume $H$ is connected. Write $|V(H)|=2k-1$. Then $\nu(H)=k-1$ and $\chi(G)=k$. Choose $v$ with $d_H(v)=\Delta(H)$ and a perfect matching of $H-v$. Define $P,C,Y,m$ as in Theorem~\ref{thm:merge}. Since $H[Y]$ has a perfect matching,
\[
\chi(G[Y])=|Y|-\nu(H[Y])=2m-m=m.
\]
Induction gives $G[Y]\succeq_s K_m$, and Theorem~\ref{thm:merge} gives $G\succeq_s K_k$.
\end{proof}

\section*{Acknowledgements}
This work was supported by the National Natural Science Foundation of China under Grant No.~12161073, the Guangdong Basic and Applied Basic Research Foundation under Grant No.~2026A1515012764, the GDUPT Talent Recruitment Project (No.~2024rcyj1007), and the Scientific Research Funds at China University of Geosciences (Wuhan) (Project No.~2026039).

\section*{Declaration of AI use}
The authors used GPT-6.1 to assist with exploring proof strategies, checking proofs, and improving the exposition. All mathematical arguments, results, and conclusions were reviewed and verified by the authors, who take full responsibility for the content of this paper.

\end{document}